\documentclass[12pt]{article}

\usepackage{amssymb,amsthm,amscd,array}

\let\ssection=\section
\renewcommand{\section}{\setcounter{equation}{0}\ssection}

\newcommand{\cA}{{\cal A}}
\newcommand{\cF}{{\cal F}}
\newcommand{\cN}{{\cal N}}
\newcommand{\bbR}{\mathbb{R}}
\newcommand{\bbC}{\mathbb{C}}
\newcommand{\bbZ}{\mathbb{Z}}
\newcommand{\ad}{\mathrm{ad}}
\newcommand{\Gr}{\mathrm{Gr}}
\newcommand{\Id}{\mathrm{Id}}
\newcommand{\ord}{\mathrm{ord}}
\newcommand{\fg}{\mathfrak{g}}
\newcommand{\fh}{\mathfrak{h}}
\newcommand{\rI}{\mathrm{I}}
\newcommand{\rZ}{\mathrm{Z}}
\newcommand{\rF}{\mathrm{F}}
\newcommand{\rC}{\mathrm{C}}
\newcommand{\rE}{\mathrm{E}}
\newcommand{\re}{\mathrm{e}}
\newcommand{\rH}{\mathrm{H}}
\newcommand{\rM}{\mathrm{M}}

\newcommand{\Res}{\mathrm{Res}}
\newcommand{\Tr}{\mathrm{Tr}}
\newcommand{\Vect}{\mathrm{Vect}}
\newcommand{\weight}{\mathrm{weight}}
\newcommand{\PD}{\Psi\mathrm{D}}

\begin{document}

\font\rus=wncyr10 scaled 1000

\newtheorem{thm}{Theorem}[section]
\newtheorem{lem}[thm]{Lemma}
\newtheorem{cor}[thm]{Corollary}
\newtheorem{pro}[thm]{Proposition}
\newtheorem{ex}[thm]{Example}
\newtheorem{rmk}[thm]{Remark}
\newtheorem{defi}[thm]{Definition}

\def\d{\delta}
\def\g{\gamma}
\def\om{\omega}
\def\r{\rho}
\def\t{\theta}
\def\a{\alpha}
\def\s{\sigma}
\def\vfi{\varphi}
\def\l{\lambda}
\def\m{\mu}

\title{Deforming the Lie algebra of vector fields on $S^1$\\
inside the Lie algebra of \\
pseudodifferential symbols on $S^1$}

\author{V. Ovsienko {\small and}
C. Roger}

\date{}

\maketitle

\date{
{\hfill{\rus Dmitriyu Borisovichu Fuksu}}

{\hfill{\rus po sluchayu shestidesyatiletiya}}
}


\thispagestyle{empty}


\begin{abstract}
We classify nontrivial deformations of the standard embedding of the
Lie algebra $\Vect(S^1)$ of smooth vector fields on the circle, 
into the Lie algebra~$\PD(S^1)$ of pseudodifferential symbols on $S^1$.
This approach leads to deformations of the central charge induced 
on $\Vect(S^1)$ by the canonical central extension of $\PD(S^1)$. 
As a result we obtain a quantized version of the second Bernoulli polynomial.
\end{abstract}

\section{Introduction}

The classical deformation theory usually deals with formal deformations of
associative (and Lie) algebras. Another part of this theory which studies
deformations of Lie algebra homomorphisms, is less known. It should be stressed,
however, that, sometimes, this second view point is more interesting and leads to
richer results.  The Lie algebra of vector fields on the circle $\Vect(S^1)$ gives
an example when such situation occurs. It is well-known that $\Vect(S^1)$ itself is
rigid, but it has many interesting embeddings to other remarkable Lie algebras that
can be nontrivially deformed. 

The aim of this article is to study the embeddings of $\Vect(S^1)$ into the Lie
algebra of pseudodifferential symbols~$\PD(S^1)$. We will classify the deformations
of the standard embedding which polynomially depend on the parameters of
deformation. It turns out that there exists a three-parameter family of nontrivial
deformations. We compute a universal explicit formula describing this family.

This work can be considered as the second part of \cite{ORP}
where deformations of~$\Vect(S^1)$ inside
the Poisson Lie algebra of Laurent series on $T^*S^1$ has been
classified. The latter Lie algebra can be considered as the semi-classical limit of
$\PD(S^1)$, so the present article is, in some sense, the quantum version of
\cite{ORP}.

It turns out that the space of infinitesimal deformations of~$\Vect(S^1)$
inside~$\PD(S^1)$ is four-dimensional. The integrability condition, i.e., the
condition of existence of a polynomial deformation corresponding to a given
infinitesimal one, distinguishes an interesting algebraic surface in the space of
parameters.

The Lie algebra of pseudodifferential
symbols has a well-known nontrivial central extension $\widehat{\PD(S^1)}$.
The restriction of this central extension to the
subalgebra~$\Vect(S^1)\hookrightarrow\PD(S^1)$ defines the famous central extension
of
$\Vect(S^1)$ called the Virasoro algebra. As an application of our results, we
obtain a ``deformed'' expression for the Virasoro central charge induced by
the deformations of the standard embedding we have constructed. It is given by a
quantized version of the second Bernoulli polynomial.

Similar results hold for the semi-direct product
$\Vect(S^1)\!\ltimes\! C^\infty(S^1)$ by the space of smooth functions on~$S^1$.

Let us also mention that a few examples of non-standard embeddings of $\Vect(S^1)$
to the Lie algebra of differential operators have been considered in \cite{FKRW}.

\section{Definitions and notations}

Let us first recall the definition of the Poisson algebra of Laurent series on $S^1$
(cf.~\cite{ORP}). Consider the space ${\bbC}[\xi,\xi^{-1}]]$ of (formal) Laurent
series of finite order in $\xi$:
$$
a(\xi)=\sum_{k\in\bbZ}a_k\xi^k,
$$
with $a_k=0$ for sufficiently large $n$. We put
\begin{equation}
{\cal A}(S^1):=C^{\infty}(S^1)\otimes{\bbC}[\xi,\xi^{-1}]]
\end{equation}
with its natural multiplication and the Poisson bracket defined as follows. Any 
element of ${\cal A}(S^1)$ can be written in the following form:
\begin{equation}
F(x,\xi)=\sum_{k\in\bbZ}f_k(x)\xi^k
\label{element}
\end{equation}
where $f_k(x)\in C^\infty(S^1)$, then for $F,G\in\cA(S^1)$
$$
\{F,G\}=
\frac{\partial F}{\partial \xi}\frac{\partial G}{\partial x}-
\frac{\partial F}{\partial x}\frac{\partial G}{\partial \xi}.
$$

The associative algebra of pseudodifferential symbols $\PD(S^1)$ has the same
underlying vector space, but the multiplication is now defined by the following
formula
\begin{equation}
(F\!\circ\! G)(x,\xi)
=
\sum_{k=0}^{\infty}\frac{1}{k!}
:
\frac{\partial^k F}{\partial \xi^k}(x,\xi)
\frac{\partial^k G}{\partial x^k}(x,\xi)
:
\label{compos}
\end{equation}
Here the double comma indicates the ordered product:
$$
:f(x)\xi^kg(x)\xi^\ell:
\; =\;
f(x)g(x)\xi^{k+\ell}.
$$

\medskip
\noindent
\textbf{Remark}.
This is not, strictly speaking, the ordered product of Wick, but rather a natural
generalization of it. Note that, here, variables $x$ and $\xi$ are no longer
commutative.

\medskip

As usual, one can consider this associative algebra $\PD(S^1)$ as a Lie algebra
with the commutator 
\begin{equation}
[F,G]=F\!\circ\! G-G\!\circ\! F.
\label{commutator}
\end{equation}

\begin{defi}
{\rm
The contraction of the algebra $\PD(S^1)$ will be constructed using the
linear isomorphism $\Phi_h:\PD(S^1)\to\PD(S^1)$ defined by
\begin{equation}
\Phi_h(f(x)\xi^k)=f(x)h^k{}\xi^k,
\label{h}
\end{equation}
where $h$ is some real number belonging to $]0,1]$.
One can modify the multiplication of~$\PD(S^1)$ according to the formula 
\begin{equation}
F\!\circ\!_hG=\Phi_h^{-1}(\Phi_h(F)\!\circ\!\Phi_h(G)).
\label{hproduct}
\end{equation}
Let us denote by $\PD_h(S^1)$ the algebra of pseudodifferential symbols equipped
with the multiplication $\circ_h$. It is clear that all the associative algebras
$\PD_h(S^1)$ are isomorphic to each other, but, for the commutator
$[F,G]_h:=F\!\circ\!_hG-G\!\circ\!_hF$, one has:
$$
[F,G]_h=
h\{F,G\}+hO(h)
$$
and, therefore, $\lim_{h\to0}\frac{1}{h}[F,G]_h=\{F,G\}.$
This is what we will mean by contraction: the Lie algebra $\PD(S^1)$ contracts to
the Poisson algebra $\cA(S^1)$.
}
\end{defi}

This notion of contraction was introduced by Wigner and In\"onu and further
developed by Levy-Nahas \cite{LN}.

Furthermore, the Lie algebra $\PD(S^1)$ admits an analogue of the Killing form,
known as the Adler trace (cf. \cite{A}). 
\begin{defi}
{\rm
Let $F(x,\xi)$ as given by (\ref{element})
be a pseudodifferential symbol, its residue $\Res\, F$ is just the coefficient
$f_{-1}(x)\in C^\infty(S^1)$ and the Adler trace is then defined by
\begin{equation}
\label{trace}
\Tr(F):=
\int_{S^1}\Res\, F\,dx=
\int_{S^1}f_{-1}(x)dx.
\end{equation}
One can readily check that $(F,G)\mapsto\Tr(F\!\circ\! G)$ defines a bilinear,
invariant, and nondegenerate form.
}
\end{defi}

\section{Statement of the problem}

The main purpose of this paper is to study deformations the canonical embedding
$\pi:\Vect(S^1)\to\PD(S^1)$ defined by 
\begin{equation}
\label{Canonical}
\pi(f(x)\partial=f(x)\xi,
\end{equation}
where $\partial=d/dx$, into a (one-)parameter family
of Lie algebra homomorphisms $\widetilde\pi_t:\Vect(S^1)\to\PD(S^1)[[t]]$, where
as usual, the symbol $[[t]]$ indicates formal power series in $t$. 

\subsection{Formal deformations}

Let us put
\begin{equation}
\label{tildepi}
\widetilde\pi_t=
\pi+\sum_{k=1}^\infty t^k\pi_k,
\end{equation}
where $\pi_k$ are linear maps $\pi_k:\Vect(S^1)\to\PD(S^1)$ and the condition
of homomorphism reads
\begin{equation}
\label{homomorphism}
\widetilde\pi_t([X,Y])=
[\widetilde\pi_t(X),\widetilde\pi_t(Y)],
\end{equation}
where the bracket in the right hand side is the natural (formal) extension to
$\PD(S^1)[[t]]$ of the Lie bracket of $\PD(S^1)$. 

\begin{defi}
{\rm
Two formal deformations $\widetilde\pi_t$ and $\widetilde\pi_t'$ are called
equivalent when there exists a formal inner automorphism
$\rI_t:\PD(S^1)[[t]]\to\PD(S^1)[[t]]$ of the form
\begin{equation}
\label{automorphism}
\rI_t=
\exp(t\,\ad F_1 +t^2\ad F_2+\cdots)
\end{equation}
with $F_i\in\PD(S^1)$, such that 
\begin{equation}
\label{equiv}
\widetilde\pi_t'=\rI_t\!\circ\!\widetilde\pi_t
\end{equation}
}
\end{defi}

\subsection{Polynomial deformations}

In many cases we will have to consider polynomial deformations instead of formal
ones: the formal series can sometimes be cut off, and the formal parameter be
replaced by complex coefficients. More precisely, a polynomial deformation has the
following form:
\begin{equation}
\label{polynom}
\widetilde\pi(c)=
\pi+\sum_{k\in\bbZ}\widetilde\pi_k(c)\xi^k,
\end{equation}
where each linear map $\widetilde\pi_k(c):\Vect(S^1)\to C^\infty(S^1)$ is polynomial
in $c\in\bbC^n$ and satisfy the following two conditions: $\widetilde\pi_k(0)=0$
and $\widetilde\pi_k\equiv0$ for sufficiently large $k$. It defines a Lie algebra
homomorphism $\widetilde\pi(c):\Vect(S^1)\to\PD(S^1)$.

\begin{defi}
{\rm
To define the notion of equivalence in the case of polynomial deformations, one
simply replaces the formal automorphism $\rI_t$ in (\ref{equiv}) by an automorphism
$\rI(c):\PD(S^1)\to\PD(S^1)$ depending on $c\in\bbC^n$ in the following way:
\begin{equation}
\label{automorphism}
\rI(c)=
\exp\left(
\sum_{i=1}^{n}c_i\,\ad F_i+\sum_{i,j=1}^{n}c_ic_j\,\ad F_{ij}+\cdots
\right),
\end{equation}
with $F_i,F_{ij},\ldots,F_{i_1,\ldots, i_k}$ belonging to $\PD(S^1)$.
}
\end{defi}

\subsection{Contraction of deformations}

Given a deformation of one of the types, (\ref{tildepi}) or (\ref{polynom}),
one can modify it so that it becomes a deformation with values in the deformed
algebra $\PD_h(S^1)$. Namely, if one
sets $\pi^h_t=\Phi_h^{-1}\!\circ\!\widetilde\pi_t$, then one has
$$
\begin{array}{rcl}
\pi^h_t([X,Y]) &=&
\displaystyle
\Phi_h^{-1}([\widetilde\pi_t(X),\widetilde\pi_t(Y)])\\[6pt]
&=&
\displaystyle
\Phi_h^{-1}([\Phi_h\pi_t^h(X),\Phi_h\pi_t^h(Y)])\\[6pt]
&=&
\displaystyle
[\pi_t^h(X),\pi_t^h(Y)]_h,
\end{array}
$$
and so $\pi^h_t:\Vect(S^1)\to\PD_h(S^1)$ is a (formal) deformation of the standard
embedding. If one sets 
$$
\pi_t=\lim_{h\to0}\pi_t^h,
$$
one obtains
$\pi_t:\Vect(S^1)\to\cA(S^1)$ which is a (formal) deformation of the standard
embedding of
$\Vect(S^1)$ into $\cA(S^1)$ as we considered and classified in \cite{ORP}. 
We will say that the deformation $\widetilde\pi_t$ contracts to the
deformation~$\pi_t$.

\medskip

We will classify the deformations of the standard embedding $\Vect(S^1)$ into
$\PD(S^1)$ and show that they all contract to some deformations of the standard
embedding $\Vect(S^1)$ into $\cA(S^1)$. This enables us to interpret these
deformations as ``quantization'' of deformations inside $\cA(S^1)$. Conversely, the
deformations inside $\cA(S^1)$ are ``semi-classical limits'' of deformations inside
$\PD(S^1)$.

\goodbreak

\section{Nijenhuis - Richardson theory of deformation}

The theory of formal and polynomial deformations of Lie algebra homomorphisms was
developed for the first time by Nijenhuis and Richardson \cite{NR}; we will make
extensive use of their methods and results, and we now recall the main points of
the cohomological apparatus for clarity and comfort of the reader.

\subsection{Infinitesimal deformations and the first
cohomology}\label{infinitesimal}

Consider the relation (\ref{homomorphism}) and take the terms of lower degree. One
has 
\begin{equation}
\label{cocycle}
\pi_1([X,Y])=
[\pi_1(X),\pi(Y)]+[\pi(X),\pi_1(Y)].
\end{equation}
In the same way, the relation (\ref{equiv}) for equivalence implies
\begin{equation}
\label{samecocycle}
\pi_1'(X)=
\pi_1(X)+[F_1,\pi(X)].
\end{equation}
Now comes a cohomological interpretation: if $\pi:\fg\to\fh$ is a Lie algebra
homomorphism, then $\fh$ is naturally a $\fg$-module through $\pi$, and
$\pi_1:\fg\to\fh$ satisfying (\ref{cocycle}) is a 1-cocycle on $\fg$ with values in
$\fg$-module $\fh$, so $\pi_1\in\rZ^1(\fg;\fh)$. Similarly, the relation
(\ref{samecocycle}) means that $\pi_1'-\pi_1$ is a coboundary, so for each formal
deformation one associates a well-defined cohomology class in $\rH^1(\fg;\fh)$.

\begin{defi}
{\rm
A map $\pi+t\pi_1:\fg\to\fh$ where $\pi_1\in\rZ^1(\fg;\fh)$ is then a Lie algebra
homomorphism up to quadratic terms in $t$, it is called an infinitesimal
deformation.
}
\end{defi} 

The same formalism works in the case of polynomial deformations, one simply has to
consider $\frac{d\pi(c)}{dc}\vert_{c=0}$ and derive the formula expressing
$\pi(c)$ as a homomorphism.

\subsection{Obstructions}

The problem is now to find higher order prolongations of these infinitesimal
deformations, there exists an obstruction theory. 

One can rewrite the relation (\ref{homomorphism}) in the following way:
$$
[\widetilde\pi_t(X),\pi(Y)]+
[\pi(X),\widetilde\pi_t(Y)]-
\widetilde\pi_t([X,Y])+
\sum_{i,j>0}[\pi_i(X),\pi_j(Y)]t^{i+j}=0.
$$
The first three terms read $(d\widetilde\pi_t)(X,Y)$ where $d$ stands for the
cohomological boundary for cochains on $\fg$ with values in $\fh$. In order to
rewrite this expression in a more compact form and to understand its cohomological
nature, we will need the following

\begin{defi}
{\rm
For arbitrary linear maps $\varphi,\varphi':\fg\to\fh$, one can associate
a bilinear map $[\![\varphi,\varphi']\!]:\fg\otimes\fg\to\fh$ in
the following way\footnote{Warning: this bracket is not a Lie algebra bracket but a
graded Lie algebra one.}
\begin{equation}
\label{cup}
[\![\varphi,\varphi']\!](X,Y)=[\varphi(X),\varphi'(Y)]+[\varphi'(X),\varphi(Y)].
\end{equation}
The relation (\ref{homomorphism}) becomes now equivalent to:
\begin{equation}
\label{defrel}
d\widetilde\pi_t+
\frac{1}{2}[\![\widetilde\pi_t,\widetilde\pi_t]\!]
=0
\end{equation}
familiar in the deformation theory (see e.g. \cite{Kon}) as the deformation
relation.
}
\end{defi}

One can now develop (\ref{defrel}) according to degree in $t$, and one obtains an
infinite series of relations
\begin{equation}
\label{obsruct}
d\pi_k+
\sum_{i+j=k}[\![\pi_i,\pi_j]\!]
=0
\end{equation}
for each $k\geq1$.
These are cohomological relations for the coboundary of
$\pi_k$, viewed as a 1-cochain on $\fg$ with values in $\fh$. 

In particular, the first nontrivial relation gives $d\pi_2+[\![\pi_1,\pi_1]\!]=0$
and one gets the first obstruction to integration of an infinitesimal deformation.
Indeed, it is quite easy to check that for any two 1-cocycles
$\gamma_1$ and $\gamma_2\in\rZ^1(\fg;\fh)$ the bilinear map
$[\![\gamma_1,\gamma_2]\!]$ is a 2-cocycle. The first nontrivial
relation (\ref{obsruct}) is precisely the condition for this cocycle to be a
coboundary. Moreover, if one of the cocycles $\gamma_1$ or $\gamma_2$ is a
coboundary, the 2-cocycle $[\![\gamma_1,\gamma_2]\!]$ is a 2-coboundary. We,
therefore, naturally obtain the following

\begin{defi}
{\rm
The operation (\ref{cup}) defines a bilinear map
$\rH^1(\fg;\fh)\otimes\rH^1(\fg;\fh)\to\rH^2(\fg;\fh)$ called the cup-product, by
analogy with the case of differential forms. This cup-product
can be extended to the whole cohomology group $\rH^*(\fg;\fh)$ (see e.g.
\cite{HS}).
}
\end{defi}

All the obstructions lie in $\rH^2(\fg;\fh)$ and they are in the image of
$\rH^1(\fg;\fh)$ through the cup-product.

So, in our case, we have to compute $\rH^1(\Vect(S^1);\PD(S^1))$ and the product
classes in $\rH^2(\Vect(S^1);\PD(S^1))$; our first task will, therefore, be to
describe
$\PD(S^1)$ as a $\Vect(S^1)$-module.

\section{Structure of $\PD(S^1)$ as a $\Vect(S^1)$-module}

\subsection{Filtration on $\PD(S^1)$}

The order of pseudodifferential operators defines a natural filtration on
$\PD(S^1)$. Recall, that the order is defined by 
$\ord(F)=\{\sup k\in\bbZ\,\vert f_k(x)\not\equiv0\}$ for every
$F(x,\xi)=\sum_{k\in\bbZ}f_k(x)\xi^k$. One sets
$$
\rF_n\left(\PD(S^1)\right)=\{F\in\PD(S^1)\,\vert\,\ord(F)\leq-n\},
$$
where~$n\in\bbZ$.
Thus, one has a decreasing filtration, 
\begin{equation}
\label{filtr}
\cdots\subset\rF_{n+1}\subset\rF_n\subset\cdots,
\end{equation}
compatible with multiplication,
if $F\in\rF_n$ and $G\in\rF_m$, then $F\!\circ\! G\in\rF_{n+m}$,
and~$\{F,G\}\in\rF_{n+m-1}$ (check the symbolic terms!).

This filtration makes 
$\PD(S^1)$ an associative filtered algebra, one can consider, as usual, the
associated graded algebra. Each quotient space $\rF_n/\rF_{n+1}$ is canonically
isomorphic to $C^\infty(S^1)$, any function $f\in C^\infty(S^1)$ induces a
pseudodifferential symbol $f\xi^{-n}$ which has a well-defined image in
$\rF_n/\rF_{n+1}$. The associated graded algebra is then 
$$
\displaystyle
\Gr\left(\PD(S^1)\right)=
\widetilde{\bigoplus_{n\in\bbZ}}\,\rF_n/\rF_{n+1}\,,
\qquad\hbox{where}\qquad
\widetilde{\bigoplus_{n\in\bbZ}}=
\Big(\bigoplus_{n<0}\Big)\oplus\Big(\prod_{n\geq0}\Big)
$$
but the induced multiplication is nothing but the restriction on the symbolic part
(i.e. the terms without derivatives in the formula (\ref{compos})) so, it is
commutative and one has
\begin{equation}
\Gr\left(\PD(S^1)\right)=\cA(S^1)
\label{gradisom}
\end{equation}
 as an associative
algebra. 

One can use this construction to recover the Poisson bracket on $\cA(S^1)$
from a purely algebraic version of symbol calculus. Let $\Phi:\PD(S^1)\to\cA(S^1)$
be the map which associates to each pseudodifferential symbol its symbolic term
(i.e. its principal symbol) and let $\Psi$ be an arbitrary section, that is,
$\Psi:\cA(S^1)\to\PD(S^1)$ satisfying $\Phi\!\circ\!\Psi=\Id_{\cA(S^1)}$, one
obtains the Poisson bracket on $\cA(S^1)$ through the formula
$$
\{F,G\}=
\Phi\left(\Psi(F)\!\circ\!\Psi(G)-
\Psi(G)\!\circ\!\Psi(F)\right)
$$

\subsection{Algebra $\cA(S^1)$ as a $\Vect(S^1)$-module}

The defined filtration is also a filtration of $\PD(S^1)$ as a
$\Vect(S^1)$-module (i.e. compatible with the natural action of $\Vect(S^1)$ on
$\PD(S^1)$). Indeed, if $X\in\Vect(S^1)$ and $F\in\rF_n(\PD(S^1))$, then 
$$
X\cdot F=[X,F]\in\rF_n(\PD(S^1)).
$$
One induces an action of
$\Vect(S^1)$ on the associative algebra $\Gr(\PD(S^1))=\cA(S^1)$ and a
simple computation shows that we recover the canonical action of $\Vect(S^1)$
on the Poisson algebra $\cA(S^1)$ considered in \cite{ORP}.

More explicitly, as a $\Vect(S^1)$-module, $\rF_n/\rF_{n+1}=\cF_{-n}$, where
$\cF_{-n}$ is the space of tensor densities of degree $n$ on $S^1$:
$$
\cF_{-n}=\left\{a(x)dx^n\,\vert\, a\in C^\infty(S^1)\right\}
$$
and the action of $\Vect(S^1)$ reads
\begin{equation}
\label{Lie}
f(x)\partial\cdot a(x)dx^n=
\left(
f(x)a'(x)+nf'(x)a(x)
\right)
dx^n.
\end{equation}
Note, that the expression in the right hand side is just the standard Lie
derivative of a tensor density along a vector field. So, finally,
\begin{equation}
\cA(S^1)=\widetilde{\bigoplus_{n\in\bbZ}}\,\cF_{-n}
\label{tensisom}
\end{equation}
as a $\Vect(S^1)$-module (cf. \cite{ORP}, Lemma 4.1).

The cohomology of $\Vect(S^1)$ with coefficients in $\cA(S^1)$ then follows from
determination of $\rH^*(\Vect(S^1);\cF_{-n})$ by D.B.~Fuchs \cite{F}, we shall
recall later the necessary results. The important fact here is that one can deduce
the cohomology for the filtered module from the cohomology of the associated graded
module, as we shall see now.

\section{Cohomology of $\Vect(S^1)$ with coefficients in the
module of pseudodifferential operators}

In this section we compute the first cohomology group $\rH^1(\Vect(S^1);\PD(S^1))$.
The Nijenhuis - Richardson theory implies (cf. Section \ref{infinitesimal}) that
this is equivalent to classification of infinitesimal deformations of the embedding
(\ref{Canonical}).

\subsection{The spectral sequence for a filtered module over a Lie algebra}

Let $\fg$ be a Lie algebra and $\rM$ a filtered module with decreasing filtration
$\{\rM_n\}_{n\in\bbZ}$ so that
$\rM_{n+1}\subset\rM_n,\,\bigcup_{n\in\bbZ}\rM_n=\rM$ and
$\fg\cdot\rM_n\subset\rM_n$. Let
$$
\Gr(\rM)=\widetilde\bigoplus_{n\in\bbZ}\,\rM_n/\rM_{n+1}
$$
be the associated
graded $\fg$-module. One can then naturally construct a filtration on the space of
cochains by setting $\rF_n(\rC^*(\fg;\rM))=\rC^*(\fg;\rM_n)$, this filtration being
obviously compatible with the Chevalley-Eilenberg differential. So, the
corresponding spectral sequence satisfies:
\begin{equation}
\rE_0^{p,q}=F_p(\rC^{p+q}(\fg;\rM))/F_{p-1}(\rC^{p+q}(\fg;\rM))
\label{spectral}
\end{equation}
and one has the
following
\begin{pro}
Let $\fg$ be a Lie algebra and $\rM$ a graded $\fg$-module as above, then one has a
spectral sequence such that 
$$
\rE_1^{p,q}=\rH^{p+q}(\fg;\Gr^p(\rM))
$$
converging to $\rH^*(\fg;\rM)$.
It induces on $\rH^n(\fg;\rM)$ a filtration coming from all the spaces
$\rE_\infty^{p,q}$ such that $p+q=n$.
\end{pro}
\begin{proof}
This is a simple adaptation of the usual construction of the spectral sequence of a
filtered complex. The reader can refer to any classical textbook on homological
algebra, for example \cite{HS}.
\end{proof}

We will use the constructed spectral sequence to compute
$\rH^1(\Vect(S^1);\PD(S^1))$ in Section \ref{Computing} below.

\subsection{$\Vect(S^1)$-cohomology with coefficients in $\cF_{-k}$}

In our case we can use the results of D.B. Fuchs to determine the cohomology with
coefficients in the graded module, namely $\rH^k(\Vect(S^1);\Gr^p(\PD(S^1)))$ for
$k=1,2$. One has:
$$
\rH^1(\Vect(S^1);\cF_{-k})=
\left\{
\begin{array}{ll}
\bbR^2,&k=0\\
\bbR,&k=1,2\\
0,&\hbox{otherwise}
\end{array}
\right.
$$
represented by the cocycles
\begin{equation}
\label{theCocycles}
\bar c_0(f\partial)=f,
\quad
c_0(f\partial)=f',
\quad
c_1(f\partial)=f''dx,
\quad
c_2(f\partial)=f'''dx^2.
\end{equation}
Now, the isomorphisms (\ref{gradisom}) and (\ref{tensisom}) imply:
\begin{equation}
\label{theCohomology}
\rH^1(\Vect(S^1);\Gr^p(\PD(S^1)))=\bbR^4
\end{equation}

For the second cohomology one has:
$$
\rH^2(\Vect(S^1);\cF_{-k})=
\left\{
\begin{array}{ll}
\bbR^2,&k=0,1,2\\
\bbR,&k=5,7\\
0,&\hbox{otherwise}
\end{array}
\right.
$$
see \cite{ORI} for explicit formul{\ae} of generators. Therefore
$\rH^2(\Vect(S^1);\Gr^p(\PD(S^1)))$ is eight-dimensional.

\subsection{Computing $\rH^1(\Vect(S^1);\PD(S^1))$}\label{Computing}

One has to check now the behavior of the above cocycles through the successive
differentials of the spectral sequence. Cocycles $\bar c_0,c_0,c_1$ and $c_2$
belong to $\rE_1^{0,1},\rE_1^{0,1},\rE_1^{1,0}$ and~$\rE_1^{2,-1}$, respectively.
The differential $d_1$ is 
$$
d_1:\rE_1^{p,q}\to\rE_1^{p+1,q}.
$$
One can check it in the following way: consider a cocycle with values in $\cA(S^1)$,
but compute its boundary as if it was with values in $\PD(S^1)$ and keep the
symbolic part of the result. This gives a new cocycle of degree equal to the degree
of the previous one plus one, and its class will represent its image by $d_1$.
The higher order differentials 
\begin{equation}
\label{thedifferential}
d_r:\rE_r^{p,q}\to\rE_r^{p+r,q-r+1}
\end{equation}
can be constructed by iteration of this procedure, the space $\rE_r^{p+r,q-r+1}$
contains the subspace coming from $\rH^{p+q+1}(\fg;\Gr^{p+1}(\rM))$.

It is now easy to see that the cocycles $\bar c_0$ and $c_0$ will survive in the
same form, we will denote them $\t_0$ and $\t_1$ when seen as cocycles with values
in $\PD(S^1)$.

For extension of the cocycles $c_1$ and $c_2$ to the algebra $\PD(S^1)$, potential
obstructions may come from $\rH^2(\Vect(S^1);\cF_\l)$ for $\l\leq7$, so we must
check the terms
$\rE_r^{1,0}\to\rE_r^{1+r,1-r}$ for $r\leq6$ and $\rE_r^{2,-1}\to\rE_r^{2+r,-r}$
for $r\leq5$. A lengthy but straightforward calculations show that both classes
survive in $\rH^1(\Vect(S^1);\PD(S^1))$; explicitly, every potential obstruction
turns out to be trivial in cohomology, so one can always find the supplementary
higher order terms to obtain two cocycles $\t_2$ and $\t_3$ coming from $c_1$ and
$c_2$ respectively. So, we have obtained
\begin{pro}
The first cohomology group $\rH^1(\Vect(S^1);\PD(S^1))$ is four-dimen\-sional and
generated by the classes of cocycles $\t_0,\t_1,\t_2$ and $\t_3$ defined above.
\label{Hpro}
\end{pro}
Another way to obtain this result is to find explicit expressions for the
nontrivial cocycles, since the cohomology group $\rH^1(\Vect(S^1);\PD(S^1))$ is
obviously upper-bounded by $\rH^1(\Vect(S^1);\cA(S^1))$.

\subsection{Explicit formul{\ae} for the $\Vect(S^1)$-cocycles}

We have already seen, the cocycles $t_0$ and $t_1$ are given by the same
formul{\ae} as $\bar c_0$ and $c_0$, namely,
\begin{equation}
\t_0(f\partial)=f,
\quad
\t_1(f\partial)=f'.
\label{t0t1}
\end{equation}
One can find the explicit formul{\ae} for $\t_2$ and $\t_3$~
\begin{eqnarray}
\t_2(f(x)\partial)
&=&
\displaystyle
\sum_{n=2}^{\infty}(-1)^{n-1}\,\frac{2(n-3)}{n}\,f^{(n)}(x)\xi^{-n+1},\\[10pt]
\displaystyle
\t_3(f(x)\partial)
&=&\sum_{n=2}^{\infty}(-1)^n\,\frac{3(n-1)}{n+1}\,f^{(n+1)}(x)\xi^{-n}.
\label{t2t3}
\end{eqnarray}

\begin{rmk}
{\rm
The symbolic terms of $\t_1(f\partial)$ and $\t_2(f\partial)$ are
$f''\xi^{-1}$ and $f'''\xi^{-2}$ respectively. Recall, that
$f\partial\mapsto f'''\xi^{-2}$ is the Souriau cocycle of the Virasoro algebra
since $\cF_2$ is the regular dual of $\Vect(S^1)$
(seeÊ\cite{K} and also \cite{I}). This class can be interpreted as infinitesimal
version of the well-known Schwarzian derivative.
}
\end{rmk}

\subsection{Relations with outer derivations of the Lie algebra $\PD(S^1)$}

The inclusion map $\Vect(S^1)\hookrightarrow\PD(S^1)$ induces the following map in
cohomology: $\rH^1(\PD(S^1);\PD(S^1))\to\rH^1(\Vect(S^1);\PD(S^1))$. The space
$\rH^1(\PD(S^1);\PD(S^1))$ of outer derivations of the Lie algebra $\PD(S^1)$ has
been determined by O.~Kravchenko and B.~Khesin \cite{KK}: it is two-dimensional and
generated by the linear operators on $\PD(S^1)$ denoted by $\ad(x)$ and
$\ad(\log\xi)$. They are defined by
\begin{eqnarray}
\displaystyle
[\ad(x),f(x)\xi^p]
&:=&
-pf(x)\xi^{p-1},\label{logx}\\[8pt]
\displaystyle
[\ad(\log\xi),f(x)\xi^p]
&:=&\sum_{n=1}^{\infty}\frac{(-1)^{n+1}}{n}\,f^{(n)}(x)\xi^{-n+p}.
\label{logxi}
\end{eqnarray}
The restriction to $\Vect(S^1)$ gives
$$
\begin{array}{rcl}
[\ad(x),f(x)\xi]
&=&
-f(x),\\[8pt]
[\ad(\log\xi),f(x)\xi]
&=&
\displaystyle
f'(x)-\frac{f''(x)}{2}\xi^{-1}+\frac{f'''(x)}{3}\xi^{-2}+\cdots\\[8pt]
&=&
\displaystyle
\left(\t_1-\frac{\t_2}{2}+\frac{\t_3}{3}\right)(f(x)\partial).
\end{array}
$$
So, the image of $\rH^1(\PD(S^1);\PD(S^1))$ in $\rH^1(\Vect(S^1);\PD(S^1))$ is
two-dimensional.

\section{Integrability of infinitesimal deformations}

Proposition \ref{Hpro} gives us all the infinitesimal deformations of the standard
embedding of $\Vect(S^1)$ into $\PD(S^1)$ given by (\ref{Canonical}). In this
section we will calculate the condition under which an infinitesimal deformation
corresponds to a polynomial one.

\subsection{Nontrivial deformation generated by cocycles $\t_0$ and~$\t_1$}

Following the Richardson-Nijenhuis theory, one has to determine the
cup-products $[\![\t_i,\t_j]\!]$ for $i,j=0,1,2,3$ of the cocycles
(\ref{t0t1}-\ref{t2t3}) in $\rH^2(\Vect(S^1);\PD(S^1))$.

First of all, since zero-order operators commute in $\PD(S^1))$, it is evident that
the the cup-products $[\![\t_0,\t_0]\!],[\![\t_0,\t_1]\!]$ and $[\![\t_1,\t_1]\!]$
vanish identically (and not only in cohomology). Therefore, the map 
\begin{equation}
f(x)\partial\mapsto f(x)\xi+\nu f(x)+\l f'(x)
\label{degone}
\end{equation}
is, indeed, a nontrivial
deformation of the standard embedding. This deformation is, in fact, both formal and
polynomial, since it is of degree one.

Real difficulties begin when we deal with integrability of the infinitesimal
deformations corresponding to the cocycles $\t_1,\t_2$ and $\t_3$ into polynomial
or formal ones.

\subsection{The integrability condition}

Consider an infinitesimal deformation of the standard embedding of
$\Vect(S^1)$ into $\PD(S^1)$ defined by the cocycles $\t_0,\t_1,\t_2,\t_3$ and
depending on three complex parameters~$c_0,c_1,c_2,c_3$:
\begin{equation}
\widetilde{\pi}(c)(f\partial)=
f\xi+c_0\t_0(f\partial)+c_1\t_1(f\partial)+c_2\t_2(f\partial)+c_3\t_3(f\partial).
\label{Infinites}
\end{equation}

\begin{thm}
\label{ConditionThm}
The infinitesimal deformation (\ref{Infinites}) corresponds to a polynomial
deformation, if and only if the following (quartic)
relation is satisfied:
\begin{equation}
\begin{array}{l}
6c_1^3c_3-3c_1^2c_2^2-18c_1c_2c_3+8c_2^3+9c_3^2+3c_1^3c_2-6c_1c_2^2\\[6pt]
\qquad\qquad
-9c_1^2c_3+18c_2c_3+2c_1c_3-5c_1^2c_2+8c_2^2+2c_1c_2=0
\end{array}
\label{quartic}
\end{equation}
\end{thm}

\noindent
Let us first prove that the condition (\ref{quartic}) is necessary for
integrability of infinitesimal deformations; we will show in Section
\ref{suffic} that this condition is, sufficient by exhibiting explicit
deformations.

\begin{proof}
The infinitesimal deformation (\ref{Infinites}) is clearly of the form:
\begin{equation}
\widetilde{\pi}(c)(f\partial)=
f\xi+c_0f+c_1f'+c_2f''\xi^{-1}+c_3f'''\xi^{-2}+\cdots,
\label{smallTerms}
\end{equation}
where ``$\cdots$'' means the terms in $\xi^{-3},\xi^{-4}$.
To compute the obstructions for integration of the infinitesimal
deformation~(\ref{Infinites}), one has to add the first nontrivial terms and impose
the homomorphism condition. So put
$$
\overline{\pi}(c)(f\partial)=
\pi(c)(f\partial)+P_4(c)f^{(IV)}\xi^{-3}+P_5(c)f^{(V)}\xi^{-4},
$$
where $P_4(c)$ and $P_5(c)$ are some polynomials in $c=(c_0,c_1,c_2,c_3)$ and
compute the difference 
$[\overline{\pi}(c)(f\partial),\overline{\pi}(c)(g\partial)]-
\overline{\pi}(c)([f\partial,g\partial])$. Collecting the terms
in $\xi^{-3}$ and $\xi^{-4}$ yields:
\begin{eqnarray}
2P_4(c)&=&2c_1c_3-c_2^2+c_1c_2-3c_3-c_2,
\label{P4}\\[4pt]
5P_5(c)&=&-c_2c_3+3c_1P_4(c)+c_2^2-6P_4(c)-c_1c_2+c_2.
\label{P5}
\end{eqnarray}
Note that these expressions do not contain $c_0$.

Let us go one step further, expand our deformation up to $\xi^{-6}$, that is, put
$\overline{\overline{\pi}}(c)(f\partial)=
\overline{\pi}(c)(f\partial)+P_6(c)f^{(VI)}\xi^{-5}$, the homomorphism condition
leads to a nontrivial relation for the parameters:
$$
\begin{array}{l}
5(-c_2P_4+4c_1P_5+c_2c_3-10P_5-c_2^2+c_1c_2-c_2)\\[6pt]
\qquad\qquad =
9\,(-2c_3^2+3c_2P_4+3c_2c_3-6c_1P_4-4c_1^2c_3+10P_4+5c_3)
\end{array}
$$
Substituting the expressions (\ref{P4}) and (\ref{P5}) for $P_4$ and $P_4$, one gets
the for\-mu\-la~(\ref{quartic}). We have thus shown that this condition
is necessary for integrability of infinitesimal deformations.~\end{proof}

\begin{rmk}
{\rm
In a cohomological interpretation, the obstructions to integrability of an
infinitesimal deformation (\ref{Infinites}) which does not satisfy the condition
(\ref{quartic}), corresponds to a nontrivial class of $\rH^2(\Vect(S^1);\cF_5)$.
}
\end{rmk}

\subsection{Introducing the parameter $h$}

One can now modify the relation in order to get a deformation in $\PD_h(S^1)$, the
scalar $h$ then appears with different powers according to the "weight" of the
successive terms in the formula (\ref{quartic}). One gets finally:
\begin{equation}
\begin{array}{rrl}
&&6c_1^3c_3-3c_1^2c_2^2-18c_1c_2c_3+8c_2^3+9c_3^2\\[6pt]
&+&h(3c_1^3c_2-6c_1c_2^2-9c_1^2c_3+18c_2c_3)\\[6pt]
&+&h^2(2c_1c_3-5c_1^2c_2+8c_2^2)\\[6pt]
&+&h^3(2c_1c_2)=0
\end{array}
\label{hquartic}
\end{equation}
This relation is a necessary integrability condition for infinitesimal
deformations~(\ref{Infinites}) in $\PD_h(S^1)$.

\begin{rmk}
{\rm
Let us remark that by setting: $\weight(c_i)=i$ (for $i=1,2,3$) and 
$\weight(h)=1$, this polynomial is homogeneous of weight 6. Moreover, put $h=0$,
one gets the polynomial we have obtained in the semi-classical
(Poisson) case (cf.~\cite{ORP}, Theorem~5.1).
}
\end{rmk}

We will now give two natural descriptions of the surface defined by the
equation~(\ref{hquartic}) in order to unveil its algebraic nature.

\subsection{Cubic curve in the space of parameters}

Let us change the parameters $c_1,c_2$ and $c_3$ and rewrite
the expression (\ref{hquartic}) in a canonical form. One checks by
an elementary straightforward calculations the
\begin{pro}
The relation (\ref{hquartic}) is equivalent to the following cubic
\begin{equation}
Y^2 = X^3 + \frac{h^2}{\!4}\,X^2,
\label{newRel}
\end{equation}
where the new parameters $X$ and $Y$ are given by
\begin{eqnarray}
X &=& c_1^2 - 2c_2 - hc_1,
\label{X}\\[4pt]
Y &=& 
\textstyle
c_1^3-3(c_1c_2-c_3) - \frac{3}{2}h(c_1^2-2c_2) + \frac{1}{2}h^2c_1,
\label{Y}
\end{eqnarray}
\end{pro}

One thus
obtains a cubic curve on the plane $(X,Y)$. This curve contracts for $h\to0$ to the
semi-cubic parabola found in \cite{ORP}.

\subsection{Rational parameterization}

There exists another, parametric, way to define the surface (\ref{hquartic}).
\begin{pro}
\label{zamenaLem}
(i) For all $\l$ and $\m\in\bbC$, the constants
\begin{eqnarray}
c_1 &=& \l+\m 
\label{zamena1}
\\[4pt]
c_2 &=& 
\displaystyle
\l\m +\frac{\l(\l-h)}{2}
\label{zamena2}
\\[6pt]
c_3 &=&
\displaystyle
\frac{\l\m(\l-h)}{2}+\frac{\l(\l-h)(\l-2h)}{6}
\label{zamena3}
\end{eqnarray}
satisfy the relation (\ref{hquartic}).

(ii) Any triple $c_1,c_2,c_3\in\bbC$ satisfying the relation (\ref{hquartic}) is of
the form (\ref{zamena1}-\ref{zamena3}) for some $\l,\m\in\bbC$.
\end{pro}
\begin{proof}
(i) The first statement can be easily checked directly. 

\medskip

(ii) To prove the converse statement, fix $c_1$ and $c_2$ as in
(\ref{zamena1}) and (\ref{zamena2}) respectively and solve (\ref{hquartic}) as a
quadratic equation with $c_3$ indeterminate. Then verify, that one of the
solutions of this equation coincides with (\ref{zamena3}) and the other one is as
follows:
$$
c_3 = \frac{\l\m(\l-h)}{2}+\frac{\l(\l-h)(\l-2h)}{6}
-\frac{\m(\m-h)(2\m-h)}{3}.
$$
But, this second solution also coincides with (\ref{zamena3}) after the involution
\begin{equation}
(\l,\m)\mapsto(\l+2\m-h,-\m+h)
\label{involution}
\end{equation}
which preserves the expressions (\ref{zamena1}) and (\ref{zamena2}).
\end{proof}

\begin{rmk}
{\rm
We have shown that $\l$ and $\m$ parameterize the surface (\ref{hquartic}); note
finally, that $X$ and $Y$ defined by (\ref{X}) and (\ref{Y}) are of the form:
$$
\begin{array}{rcl}
X &=& \m^2-h\m,\\[6pt]
Y &=& \m^3-\frac{3}{2}h\m^2+\frac{1}{2}h^2\m,
\end{array}
$$
and so $\m$ is a parameter on the curve (\ref{newRel}).
}
\end{rmk}

Now, the Richardson-Nijenhuis theory prescribes us to compute the second cohomology
group $\rH^2(\Vect(S^1);\PD(S^1))$ in order to obtain the complete information
concerning the cohomological obstructions. This, however, seems to be a quite
difficult problem. If one tries to compute $\rH^2(\Vect(S^1);\PD(S^1))$ using the
spectral sequence, just as we did it in Section \ref{Computing} for the first
cohomology group, one has to check $\rH^3(\Vect(S^1);\cF_n)$ for $n\leq15$ (see
\cite{F}, p.176) in order to decide whether the nontrivial classes from
$\rH^2(\Vect(S^1);\cA(S^1))$ will survive or not. We shall not do that;
an explicit construction of polynomial deformations will
allow us to avoid the standard framework.

\section{The universal formula for nontrivial deformations}\label{suffic}

In this section we will give explicit formul{\ae} of genuine polynomial
deformations of the standard embedding (\ref{Canonical}) for all infinitesimal
deformations satisfying the condition (\ref{hquartic}). According to the
Richardson-Nijenhuis theory, any polynomial deformation is equivalent to a
deformation from this class.

\subsection{Construction of deformations}\label{Devia}

Outer derivations (\ref{logx}) and (\ref{logxi}) can be integrated in one-parameter
families of outer automorphisms denoted by $\Phi_\nu$ and $\Psi_\m$ respectively,
and defined by
\begin{eqnarray}
\Phi_\nu(F)&=&\re^{x\nu}\circ F\circ\re^{-x\nu}\label{auto1}
\\[8pt]
\Psi_\m(F)&=&\xi^\m\circ F\circ\xi^{-\m}.
\label{auto}
\end{eqnarray}
Note that these formul{\ae} should be understood as Laurent series in $\xi$
(depending on the parameters $\nu$ and $\m$ respectively); one computes the product
in (\ref{auto1}) and~(\ref{auto}) formally, using the composition formula
(\ref{compos}), as if $\re^{x\nu}$ and $\xi^\m$ were elements of $\PD(S^1)$. One
obtains, in particular: 
$$
\Phi_\nu(f\xi)=f\xi+\nu f
\qquad\hbox{and}\qquad 
\Psi_\m(f\xi)=\sum_{n=0}^{\infty}{\m\choose
n}f^{(n)}\xi^{-n+1}.
$$

Let us apply the automorphism (\ref{auto}) to the elementary
deformation (\ref{degone}):
\begin{equation}
\begin{array}{rcl}
\widetilde{\pi}_{\l,\m}(f\partial)
&=&
\Psi_\m(f\xi+\l f')\\[12pt]
&=&
\displaystyle
f\xi+(\l+\m)f'\\[8pt]
&&
\displaystyle
+\left(
\l\m +\frac{\l(\l-h)}{2}
\right)f''\xi^{-1}\\[14pt]
&&
\displaystyle
+\left(
\frac{\l\m(\l-h)}{2}+\frac{\l(\l-h)(\l-2h)}{6}
\right)f'''\xi^{-2}\\[12pt]
&&
+\cdots
\end{array}
\label{TheUniversDef1}
\end{equation}
Since $\Psi_\m$ is an automorphism, it is, indeed, a polynomial
deformation of the embedding (\ref{Canonical}) for any $\l$ and $\m\in\bbC$.
In the same way, one also obtains a more general class of polynomial
deformations:
\begin{equation}
\widetilde{\pi}_{\l,\m,\nu}(f\partial)
=
\Psi_\m\left(f\xi+\nu f+\l f'\right)
=
(\Psi_\m\circ\Phi_\nu)\left(f\xi+\l f'\right).
\label{TheUniversDef}
\end{equation}

\begin{lem}
Every infinitesimal deformation (\ref{Infinites}) satisfying the condition
(\ref{hquartic}) can be realized as the infinitesimal part of 
the polynomial deformation $\widetilde{\pi}_{\l,\m,\nu}$ for some
$\l,\m$ and $\nu\in\bbC$.
\end{lem}
\begin{proof}
The first coefficients in (\ref{TheUniversDef1})
coincide with (\ref{zamena1}-\ref{zamena3}). The formula (\ref{smallTerms})
shows then that the parameters $c_1,c_2,c_3$ of the infinitesimal
part of (\ref{TheUniversDef1}) (and (\ref{TheUniversDef})) are precisely given by
these expressions. Now, Proposition \ref{zamenaLem} implies that any infinitesimal
deformation (\ref{Infinites}) satisfying the condition (\ref{hquartic}) is the
infinitesimal part of one of the deformations (\ref{TheUniversDef}).
\end{proof}

\begin{rmk}
{\rm
(a) The constructed deformations are not equivalent to each other for different
values of the parameters $\l,\m$ and $\nu$, since the corresponding infinitesimal
deformations are given by non-cohomological cocycles. This is due to the fact that
$\Phi_\nu$ and $\Psi_\m$ are outer automorphisms.

(b) After the contraction procedure (for $h\to0$) applied to the
deformations (\ref{TheUniversDef1},\ref{TheUniversDef}) one gets precisely the
deformations, found in the semi-classical case in \cite{ORP} (see Theorem 5.1).
}
\end{rmk}

\subsection{Proof of Theorem \ref{ConditionThm}}

We have constructed a polynomial deformation corresponding to any infinitesimal
deformation satisfying the condition (\ref{quartic}) (or (\ref{hquartic})). This
implies that the condition~(\ref{quartic},\ref{hquartic}) is not
only necessary, but also sufficient for integrability of infinitesimal
deformations and completes the proof of Theorem \ref{ConditionThm}.

\section{Variation of the Virasoro central charge}

The space $\rH^2(\PD(S^1);\bbC)$, classifying central extensions of the Lie algebra
$\PD(S^1)$ has been determined in \cite{KK} (see also \cite{R}). It is
two-dimensional, every outer derivation $\d_i$ (for $i=1,2$) given by the
formul{\ae} (\ref{logx}) and (\ref{logxi}) defines a 2-cocycle with scalar values
through the formula
\begin{equation}
c_i(F,G)=\int_{S^1}\Res(\d_i(F)\circ G).
\end{equation}

It is well known that the space $\rH^2(\Vect(S^1);\bbC)$ is one-dimensional and
generated by the famous Gelfand-Fuchs cocycle, that we will define by the formula
\begin{equation}
c(f\partial,g\partial)=\frac{1}{12}\int_{S^1}f'''g\,dx.
\label{GelFuc}
\end{equation}
\begin{rmk}
{\rm
The normalization term $1/12$ is traditional: it is very natural because of the
values of central charges one encounters among the Virasoro algebra representations.
}
\end{rmk}

We will consider the cocycle $c_1\in\rZ^2(\PD(S^1);\bbC)$, associated with
$\d_1=\log\xi$. A lengthy, but easy computation then proves the following
\begin{pro}
The restriction of the cocycle $c_1$ to $\Vect(S^1)\hookrightarrow\PD(S^1)$ with
respect to the embedding (\ref{TheUniversDef})
is $\widetilde{\pi}_{\l,\m,\nu}^*(c_1)=(-12\l^2+12\l-2)\,c,$ and in the case
of $h$-deformed algebra $\PD_h(S^1)$:
\begin{equation}
\widetilde{\pi}_{\l,\m,\nu}^*(c_1)=(-12\l^2+12h\l+4-6h)\,c,
\label{Bernoulli}
\end{equation}
\end{pro}

\begin{rmk}
{\rm
One then recovers for $h=1$ the famous (second) Bernoulli polynomial which appears
in central charges of some representations of the Virasoro algebra, as well as in
the context of moduli spaces (see \cite{ACKP}). The formula (\ref{Bernoulli}) is,
therefore, its ``quantized'' version.
}
\end{rmk}

\section{Deformations of the embedding of $\Vect(S^1)\!\ltimes\!\cN$}

Denote $\cN$ the space of $C^\infty$-functions on $S^1$. The semi-direct product
$\Vect(S^1)\!\ltimes\!\cN$ is then the Lie algebra of differential operators of order
smaller or equal to 1. One can look for deformations of the embedding of this Lie
algebra into $\PD(S^1)$.

Tedious, but without serious difficulties, computations allow to prove that every
deformation of the canonical embedding can be written in the form
\begin{equation}
\widehat{\pi}_{\l,\m,\nu}\Big(f(x)\partial+a(x)\Big)=
f(x)\xi+\nu f(x)+\l f'(x)+(\m+1)a(x).
\label{semiimbed}
\end{equation}
One can now compute the induced map on second cohomology.

As above, the space $\rH^2(\Vect(S^1)\!\ltimes\!\cN;\bbC)$ is well known (see
\cite{ACKP}) and admits three generators: the first one is the Gelfand-Fuchs cocycle
still denoted $c$, and the following two cocycles:
\begin{eqnarray}
\widetilde{c}&=&\int_{S^1}(f''b-g''a)dx,\\[6pt]
\widetilde{\widetilde{c}}&=&\int_{S^1}(a'b-ab')dx.
\end{eqnarray}
One has the following
\begin{pro}
The restriction of the cocycle $c_1$ to
$\Vect(S^1)\!\ltimes\!\cN\hookrightarrow\PD(S^1)$ with respect to the embedding
(\ref{semiimbed}) is as
follows
$$
\widehat{\pi}_{\l,\m,\nu}^*(c_1)
=(-12\l^2+12\l-2)\,c+
(\nu+1)(\l-\frac{1}{2})\,\widetilde{c}+
\frac{(\nu+1)^2}{2}\,\widetilde{\widetilde{c}}
$$
and in the case of $h$-deformed algebra $\PD_h(S^1)$:
$$
\widehat{\pi}_{\l,\m,\nu}^{h\;*}(c_1)
=(-12\l^2+12h\l+4-6h)\,c+
(\nu+1)(\l-\frac{h}{2})\,\widetilde{c}+
\frac{(\nu+1)^2}{2}\,\widetilde{\widetilde{c}}
\label{BernoulliFirst}
$$
\end{pro}

\begin{rmk}
{\rm
The coefficient $\l-\frac{1}{2}$ is also a Bernoulli polynomial, the first one.
In terms of $W$-algebras, we computed central charges of $W_1$ and $W_2$.
}
\end{rmk}

\bigskip

\noindent
{\bf Acknowledgments}. We are grateful to Ch. Duval and E. Ragoucy for stimulating
discussions.


$$
\begin{array}{lcl}
\hbox{Valentin OVSIENKO}
&&
\hbox{Claude ROGER}
\\
\hbox{{\small C.N.R.S.}}
&&
\hbox{{\small Institut Girard Desargues, UPRESA CNRS 5028}}
\\
\hbox{{\small  Centre de Physique Th\'eorique}}
&&
\hbox{{\small Universit\'e Claude Bernard - Lyon I}}
\\ 
\hbox{{\small  Luminy-Case 907}}
&&
\hbox{{\small 43 bd. du 11 Novembre 1918}}\\
\hbox{{\small  F-13288 Marseille Cedex 9, France}}
&&
\hbox{{\small 69622 Villeurbanne Cedex, France}}\\
\hbox{{\small  mailto: ovsienko@cpt.univ-mrs.fr}}
&&
\hbox{{\small mailto: roger@desargues.univ-lyon1.fr}}
\end{array}
$$

\end{document}